\documentclass[review,authoryear,11pt]{elsarticle}

 \biboptions{sort&compress}
\usepackage{times}
\usepackage{float}
\usepackage[colorlinks=true,linkcolor=blue, pdfstartview=FitH, breaklinks=true]{hyperref}
\usepackage{longtable,lscape,verbatim}
\usepackage{amsmath}
\usepackage{amssymb}
\usepackage{epsfig}
\usepackage[ruled]{algorithm2e}
\allowdisplaybreaks[1] 
\usepackage{natbib}
\usepackage{graphicx}
\usepackage{booktabs}
\usepackage{tabularx}
\usepackage{lscape}
\usepackage{url}
\usepackage[ruled]{algorithm2e}
\usepackage{booktabs}
\usepackage{color}
\usepackage{url}
\usepackage[ruled]{algorithm2e}
\usepackage{booktabs}
\usepackage{color}
\usepackage{stmaryrd}
\usepackage{amsmath}
\usepackage{setspace}

\usepackage{natbib}
\newtheorem{thm}{Theorem}[section]

\newtheorem{cor}[thm]{Corollary}

\newcommand{\Real}{\mathop{\rm I\kern-.2emR}}
\newcommand{\set}[1]{\left\{{#1}\right\}}

\begin{document}

\title{
Recovery-to-Efficiency: A New Robustness Concept\\for Multi-objective Optimization under Uncertainty
}


%

\author[inria,univlille]{El-Ghazali Talbi}
\ead{el-ghazali.talbi@univ-lille1.fr}
\author[LAMIH]{Raca Todosijevi\'c \corref{mycorrespondingauthor}}
\ead{racatodosijevic@gmail.com}

%
%
%
\address[LAMIH]{LAMIH UMR CNRS 8201\& INSA Hauts De France - Universit\'e Polytechnique Hauts-de-France, 59313 Valenciennes Cedex 9, France}
\address[inria]{INRIA Lille-Nord Europe, 40 Avenue Halley 59650 Villeneuve d'Ascq, France}
\address[univlille]{Univ. Lille, CNRS, Centrale Lille, UMR 9189 – CRIStAL, F-59000 Lille, France}

\begin{abstract}
This paper presents a new robustness concept for uncertain multi-objective optimization problems.
More precisely, in the paper the so-called \emph{recovery-to-efficiency} robustness concept is proposed and investigated.
Several approaches for generating recovery-to-efficiency robust sets in the context of multi-objective optimization are proposed as well. An extensive experimental analysis is performed to disclose differences among robust sets obtained using different concepts as well
as to deduce some interesting observations. For testing purposes, instances from the bi-objective knapsack problem are considered.
\smallskip

\noindent {\bf Keywords.}  Uncertainty, Robustness measure, Recovery-to-Efficiency, Multi-objective Optimization

\end{abstract}
\maketitle

\section{Introduction}\label{intro}

A multi-objective optimization problem in its general form can be formalized as follows:
\begin{align}
\min~			&	f(x)		\tag{$P$}	\\
\text{subject to~}	&	x\in X	\notag
\end{align}
where $X$ 
is the \emph{feasible set in decision space}, and $f=(f_1,f_2, \dots, f_n):X\rightarrow \Real^n$ is a function vector to be ``minimized" subject to some constraints.
%
A solution $x \in X$ is dominated by a solution $x^\prime \in X$, denoted as $x \prec x^\prime$,
if $f_i(x^\prime) \leqslant f_i(x) $ for all $i \in \set{1,\ldots,n}$, with at least one strict inequality.
A solution $x \in X$ is \emph{efficient} if there does not exist any other solution $x^\prime \in X$ such that~$x \prec x^\prime$.
The set of all efficient solutions is the \emph{efficient set}. 
One of the most challenging issues in multi-objective optimization is to identify the efficient set, or a good approximation of it for large-size instances.

The above problem is usually referred to as a \emph{deterministic} multi-objective optimization problem, because all parameters are assumed to be known with certainty.
However, in real-world applications, the objective functions may rely on parameters which are uncertain, and for which is only known that stem from some uncertainty set. For example, in multi-objective shortest paths, while the length of paths is certain and known exactly,  the travel time along a path is often unknown beforehand, and usually depends on many factors. Hence, in real-world applications, it is preferable to consider the uncertain multi-objective shortest path problem rather than the deterministic one. In order to include such uncertainties, instead of the (deterministic) problem (P) it is hence desirable to consider the following parameterized family of problems:
\begin{align}
\min~			&	f(x,\xi)	\tag{$P(\xi)$}	\\
\text{subject to~}	&	x\in X	\notag
\end{align}
%
where $f(\cdot,\xi)=(f_1(\cdot,\xi),f_2(\cdot,\xi), \dots, f_n(\cdot,\xi)):X\rightarrow \Real^n$.
The (unknown) parameter $\xi$ actually represents a scenario that may occur. It is assumed that  $\xi$  varies within
a  given  uncertainty  set ${\cal U}\subseteq \Real^n$. Hence, the  uncertain  optimization  problem  corresponding  to
$P$ is denoted as $(P(\xi):\xi \in {\cal U})$.

The above uncertain  optimization  problems are the topics of the \emph{robust optimization} field. While robust optimization problems
and uncertainly-handling optimization concepts are very well studied in the single-objective case,
multi-objective robust optimization problems are still an under-explored area. The first concept to deal with the uncertainty in single-objective robust optimization  is  called minmax (or strict)
robustness, originally introduced by \citet{soyster1973technical}. According to this concept, a solution is called robust if it
is feasible in each scenario and if it minimizes the \emph{worst}
objective function value with respect to all scenarios.  This concept has been later revisited and extensively studied in, e.g., \citet{ben1998robust,ben1999robust,ben2000robust,ben2009robust}.
Similar concept to minmax robustness is the absolute (relative) regret~\citep{kouvelisrobust}, 
where the feasibility of a robust solution with respect to all scenarios is, once again, required, but the goal is now to minimize the absolute (relative) deviation from the optimal solution in all scenarios.
However, all these concepts are often argued 
as too conservative, since a robust solution is
required to be feasible for each scenario. This observation inspired researchers to develop alternative, less conservative, concepts. 
Those include
light robustness~\citep{Fischetti2009} and its generalizations~\citep{schobel2014generalized},
such as adjustable robustness~\citep{ben2004adjustable},
soft robustness~\citep{ben2010soft},
or recovery robustness~\citep{Liebchen2009,erera2009robust,goerigk2014recovery}.
The reader is referred to \citet{goerigk2015algorithm,gabrel2014recent} for an overview of these (single-objective) robustness concepts.

Nevertheless, the theory of robustness dealing with multi-objective optimization is rather poor, although uncertainty appears in many multi-objective problems.
The research attempts in this line 
mainly aim at bringing the concepts originally proposed for single-objective robust optimization
into the multi-objective case. For instance, inspired by minmax robustness in single-objective optimization,
\citet{kuroiwa2012robust} replace the objective functions by their respective worst cases over all scenarios, and
hence obtain a deterministic multi-objective optimization problem whose efficient solutions are
so-called robust.
Similarly,  by extending the concept of point-based minmax robust efficiency for the single-objective
case, the concepts of set-based minmax robust efficiency \citep{ehrgott2014minmax} and Hull-based minmax robust efficiency \citep{BOKRANTZ2017682} are obtained.
More recently, 
\citet{Ide2016} proposed a light robustness concept for multi-objective robust optimization by extending the light robustness concept from the single-objective case. \cite{RAITH2018628} extended the concept of \citet{bertsimas2003robust} to the multi-objective setting, while \citet{botte2019dominance} introduced multi-scenario efficiency for robust multi-objective programming.
\cite{ide2014concepts} and \cite{ide2014relationship} elaborated the relation of uncertain multi-objective  optimization problems to the field of set-valued optimization.

Uncertain multi-objective linear  problems have been studied in several papers in terms of optimality and duality conditions,  robust counterparts and numerically tractable optimality conditions. These papers include \cite{Wang2018, dranichak2018highly,doolittle2016note,goberna2015robust,HASSANZADEH2013357,goberna2014robust}. On the other hand,
uncertain multi-objective convex  problems have been studied in \citep{GOBERNA201840,kuroiwa2014robust}. \citet{lee2018optimality} studied optimality conditions and duality theorems for robust semi-infinite multi-objective optimization problems, while \citet{chuong2016optimality, lee2015nonsmooth,zamani2015robustness} did so for robust nonsmooth/nonconvex  multi-objective optimization problems.

The robust multi-objective optimization has many real-world applications emerging in engineering, business, and management. Some of them are: robust internet routing \citep{doolittle2012robust}, forestry management \citep{palma2010bi}, portfolio management \citep{fliege2014robust, xidonas2017arobust, xidonas2017brobust, fakhar2018nonsmooth}, a multi-objective robust optimisation of carbon and glass fibre-reinforced hybrid composites under flexural loading \citep{kalantari2016multi}, aircraft route guidance and hazardous materials routing problems \citep{kuhn2016bi}, game theory \citep{yu2013robust},  and wood industry \citep{ide2015application}.
For an 
survey on existing robustness concepts for multi-objective robust optimization, we refer the reader to \citet{Ide2016}. 


In this paper, we follow the research line of extending a robustness concept from single-objective optimization to multi-objective optimization.
More precisely, we propose and investigate the so-called \emph{recovery-to-efficiency} robustness concept, which is an extension and generalization of recovery-to-optimality robustness from \citet{goerigk2014recovery} in the single-objective setting.
Hence, the main purpose of this work is to attempt further reducing the gap between the fields of robust and multi-objective optimization.
In this work, we propose several approaches for generating recovery-to-efficiency robust sets in the context of multi-objective optimization. The difference among the robust sets obtained using different approaches is revealed trough a detailed experimental analysis. For testing purposes, instances from the bi-objective knapsack problem are considered. In brief, the contributions of the paper may be summarized as follows:
\begin{itemize}
\item	New robustness concept for uncertain multi-objective optimization problems is proposed;
\item	Several approaches for generating recovery-to-efficiency robust sets are elaborated;
\item	An extensive testing on instances from the bi-objective knapsack problem is performed;
\item	The difference among the robust sets obtained using different approaches is revealed.
\end{itemize}
%

The rest of the paper is organized as follows.
In Section~\ref{sec:rte}, we introduce the recovery-to-efficiency approach as a robustness concept for multi-objective optimization.
In Section~\ref{sec:robustsets}, we give some fundamental properties from recovery-to-efficiency robustness.
In Section~\ref{sec:exp}, we provide an experimental analysis of our approach by measuring its effects when solving the robust bi-objective knapsack problem.
In Section~\ref{sec:conclu}, we conclude the paper and discuss ideas for future works.

\section{Recovery-to-Efficiency}\label{sec:rte}

In this section,
we  introduce  a  new  robustness approach for uncertain multi-objective optimization problems, strictly related to the  less conservative approach from single-objective robust optimization, namely recovery-to-optimality \citep{goerigk2014recovery}. In order to define this new concept we need a \emph{distance function} $d: X \times X \rightarrow \Real^n$,
which represents the \emph{recovery cost} to transform one solution into another one. It is not required that such function satisfies
any property of a norm or a metric, but rather reflects accurately the required effort in order to modify one solution into another one. As a consequence, such function is typically specific to the problem at hand. For instance,  in timetabling,  it may present the  increase  in travel times from one timetable to another, or for
knapsack problems, it may be the Hamming distance between two solutions, i.e., the number of items to be removed or inserted to go from one solution to another.

For each scenario $\xi$, solving problem $P(\xi)$ means identifying the set of efficient solutions. Let $x(\xi)$ denote
the efficient set of problem $P(\xi)$. In order to find $x(\xi)$, many techniques have been proposed in the literature. Among them,
the widely-used scalarizing technique consists in choosing vectors $\lambda=(\lambda_1,\lambda_2,\dots \lambda_n) \in \Real^n, \sum_{i=1}^n \lambda_i=1, \lambda_i>0$ and solving a series of sub-problems
$$(P(\xi, \lambda)) \, \min F(f(x,\xi), \lambda)$$
subject to
$$x\in X$$
where $F$ stands for functional.

The most common scalarizing techniques are weighted sum and Chebychev scalar functions. The weighted sum consists of choosing vectors $\lambda=(\lambda_1,\lambda_2,\dots \lambda_n) \in R^n,  \sum_{i=1}^n \lambda_i=1, \lambda_i>0$ and solving series of problems

$$(P(\xi, \lambda)) \, \min \sum_{i=1}^n \lambda_if_i(x,\xi)$$
subject to
$$x\in X$$
in order to determine a discretization of the efficient set.

On the other hand the Chebychev sum consists of choosing vectors $\lambda=(\lambda_1,\lambda_2,\dots \lambda_n) \in R^n,  \sum_{i=1}^n \lambda_i=1, \lambda_i>0$ and solving series of sub-problems

$$(P(\xi, \lambda)) \, \min \max_{i=1,2, \dots n} \lambda_i(f_i(x,\xi)-h_i)$$
subject to
$$x\in X$$
in order to determine an efficient set approximation, where $(h_1,\dots, h_n)$ represents reference point (which may be e.g., the nadir objective vector).

Assuming that, for each given scenario $\xi$, we can compute the set of efficient solutions $x(\xi)$, we seek for a set $x$ from which each set $x(\xi)$, $\xi \in {\cal U}$, may be easily accessed, analogously to the concept presented by \citet{goerigk2014recovery}. The notion of ``easily accessible'' here refers to finding the solution set $x$ from which the decision maker can choose one solution, and after a given scenario $\xi$ becomes known with certainty, he/she can easily transform the chosen solution into an efficient solution from set $x(\xi)$. Then denote the solution set $x$ as the set of robust recoverable solutions. As already pointed out, a similar concept already exists in single-objective optimization \citep{goerigk2014recovery}, where 
the authors seek for a solution $x$ that can easily be transformed into an optimal solution for each considered scenario.
The multi-objective concept that we present here actually generalizes the one from \citet{goerigk2014recovery} for any number of objectives.

Let us now assume that each set $x(\xi)$ is generated by means of a scalarizing approach, such that $x(\xi)=\{x(\xi, \lambda):\lambda\in \Lambda\}$, where
$x(\xi, \lambda)$ refers to the solution of  problem $P(\xi, \lambda)$ and $\Lambda=\{\lambda\in \Real^n| \sum_{i=1}^n \lambda_i=1, \lambda_i>0\}$.  Then, one approach to generate the robust set $x$ of solutions is to generate point $x(\lambda)$ for each $\lambda\in \Lambda$,  which is in some sense close to the points $x(\xi, \lambda)$, $\xi \in {\cal U}$. If the closeness is measured as the largest distance, the problem of finding   $x(\lambda)$ turns to solving the following optimization problem:
\begin{align}
\min~			&	\sup_{\xi \in {\cal U}}d(x(\lambda), x(\xi, \lambda))		\tag{Rec-Eff center$(\lambda)$}	\\
\text{subject to~}	&	x(\xi, \lambda)\in x(\xi), \; \xi \in {\cal U}			\notag	\\
				&	x(\lambda) \in X								\notag
\end{align}
%
%
If the objective in $(\text{Rec-Eff center}(\lambda))$, which aims at minimizing the maximum distance, is replaced by an objective aiming at minimizing the sum
of distances (i.e., $\min~				\sum_{\xi \in {\cal U}}d(x(\lambda), x(\xi, \lambda)$) , we refer to the resulting problem as $(\text{Rec-Eff median}(\lambda))$.
It is important to notice that,
once the decision maker has chosen a solution $x(\lambda)\in x$ from the provided robust set, he/she will automatically know how to react if a certain scenario $\xi \in {\cal U}$ is realized. The proper reaction would be to change the solution $x(\lambda)$ by the solution $x(\lambda, \xi)$.

However, in the case that the problem $P(\xi, \lambda)$ has more than one optimal solution, the better option would be to determine simultaneously $x(\lambda)$ and $x(\xi, \lambda)$.
Let us denote the 
value of a solution $x(\xi, \lambda)$ with respect to the problem $P(\xi, \lambda)$ as $g(x(\xi, \lambda))$,
and the optimal value of $P(\xi, \lambda)$ as $g_{opt}(\xi, \lambda)$.
We may define a new set of problems as follows:
\begin{align}
\min~			&	\sup_{\xi \in {\cal U}}d(x(\lambda), x(\xi, \lambda))		\tag{Rec-Eff center\_opt$(\lambda)$}	\\
\text{subject to~}	&	g(x(\xi, \lambda)) \leq g_{opt}(\xi, \lambda), \;  \xi \in {\cal U}		\notag	\\
				&	x(\lambda),x(\xi, \lambda)  \in X, \;  \xi \in {\cal U}							\notag
\end{align}
%
which sought to determine, simultaneously, the robust solution $x(\lambda)$ and the efficient solution $x(\xi, \lambda)$ for each scenario. In the case that each problem $P(\xi, \lambda)$ has a unique solution, problems  $(\text{Rec-Eff center}(\lambda))$ and $(\text{Rec-Eff center\_opt}(\lambda))$ are equivalent. Otherwise, problem $(\text{Rec-Eff center\_opt}(\lambda))$ aims at producing the set of  $x(\lambda)$ with the cheapest recovery-to-efficiency costs. Notice also that equation:
$$g(x(\xi, \lambda)) \leq g_{opt}(\xi, \lambda), \;  \xi \in {\cal U}	$$
can be replaced by the following one, allowing the ``efficient solutions" of each scenario to be within some tolerance $\epsilon$ with respect to the optimal objective value:
$$g(x(\xi, \lambda)) \leq (1+\epsilon) \; g_{opt}(\xi, \lambda), \;  \xi \in {\cal U}.$$
Analogously, to problem $(\text{Rec-Eff center\_opt}(\lambda))$ $(\text{Rec-Eff median\_opt}(\lambda))$ may be defined.

\medskip
\noindent
{\underline {\bf Single objective case:}} In the case of single-objective optimization, i.e., when $n=1$, we have that solving a problem $P(\xi, \lambda)$ for each $\lambda$ yields an optimal solution for the scenario $\xi \in {\cal U}$. Hence, solving problem $(\text{Rec-Eff center}(\lambda))$ results in finding a solution $x$ that minimizes the worst case distance to the optimal solutions for all scenarios. This means that the proposed concept includes the concept presented by \citet{goerigk2014recovery} as a special case for $n=1$. However, there are some differences with respect to the approach  presented in \citet{goerigk2014recovery}. \citet{goerigk2014recovery} consider that for each scenario all optimal solutions are known and aim to determine a  ``recovery-to-optimality" solution with respect to all optimal solutions in each scenario. However, enumerating all optimal solutions of a problem sometimes is not an easy task. Hence, in this paper in problems $(\text{Rec-Eff center}(\lambda))$ and $(\text{Rec-Eff median}(\lambda))$ we assume to have one optimal solution for each scenario (calculated by certain algorithm). Further, to avoid to be misleaded by considering only one optimal solution for each scenario, we propose an alternative approach where  the robust solution  and the ``efficient solutions" are determined simultaneously through problems $(\text{Rec-Eff center\_opt}(\lambda))$ and $(\text{Rec-Eff median\_opt}(\lambda))$. Also, we introduce another dimension of flexibility through constraints  $$g(x(\xi, \lambda)) \leq (1+\epsilon) \; g_{opt}(\xi, \lambda), \;  \xi \in {\cal U},$$
which allow that  ``efficient solutions" of each scenario are not necessarily the optimal ones but high quality ones (i.e., solutions having acceptable tolerance from the optimal solution value) .

\section{Generating Robust Sets by Sampling}
\label{sec:robustsets}

Finding  a recovery-to-efficiency  robust  solution  as  an optimal  solution of $(\text{Rec-Eff center})$ (resp. $(\text{Rec-Eff median})$) theoretically requires to find a point $x$ which minimizes the maximum (resp. sum of) distance(s) to an infinite number of points. 
The distance between two points represents
the recovery costs for changing one solution to another one. However, in practical settings, this is not a convenient approach because
such process may be time consuming. Similarly, accurately finding an efficient set $x(\xi)$ requires a large (sometimes infinite) number of different weighting coefficient vectors ($\lambda$'s), that also result into high computational burden. Nevertheless, in this paper, we come out with an approach which generalizes the one
proposed by \citet{goerigk2014recovery} for the single-objective case, and which is able to generate solutions with
low recovery costs. The proposed approach actually results into heuristic which selects a subset of scenarios and chooses a limited set of weighting coefficient vectors for which the robust set $x$ is to be created.
Algorithm \ref{alg:robustset} sketches the steps of the proposed approach,
where we assume that  an algorithm ${\cal A}$ is available for solving problem $P(\xi, \lambda)$, i.e., ${\cal A}$ takes a  given problem $P(\xi, \lambda)$ as an input and returns an optimal solution as an output.

 \begin{algorithm}[h]
 \label{alg:robustset}
\caption{Generate robust set by sampling.}
 \hspace*{10pt}  \parbox{7in} {{\bf Function} {\tt Generate\_set$(P,{\cal U})$}\;
\nl Choose set of scenarios $\{\xi^1,\xi^2,\dots, \xi^m\}\subset{\cal U}$\; \label{step:scenarios}
\nl Choose set of lambda vectors $\Lambda=\{\lambda^1,\lambda^2,\dots, \lambda^k\}$\;
\nl \For{$i=1:m$}
    {
     $x(\xi^i)=\emptyset$\;
      \For{$j=1:k$}
    {
         $x(\xi^i)_j\leftarrow{\cal A}(P(\xi^i, \lambda^j))$\;
        $x(\xi^i)\leftarrow x(\xi^i) \cup \{x(\xi^i)_j\}$\;
      }
    }
\nl $x=\emptyset$\;
\nl \For{$j=1:k$}
    {
      \nl Calculate $x^\star\in X$ minimizing the  recovery costs \\to the points $\{x(\xi^1)_j, \dots x(\xi^m)_j\}$\; \label{step:x*}
      \nl $x \leftarrow x\cup \{x^\star\}$

    }
\nl \Return x\;
}

\end{algorithm}

Similar to \citet{goerigk2014recovery}, under specific conditions, the optimal solution $x^\star\in X$ in Step \ref{step:x*} may be determined as the center of a finite number of points, which further allows a reduction of the number of  scenarios to be taken into account. These properties are stated in the next two theorems for problem $(\text{Rec-Eff center})$ and their proofs may be deduced analogously to the proofs provided in \citet{goerigk2014recovery}.

\medskip

\begin{thm}
 Let $(P(\xi:\xi \in {\cal U}))$ be an uncertain multi-objective optimization problem with problem instances $P(\xi,\lambda_j)$ which
admit  a  unique  optimal  solution for every $\xi \in {\cal U}$ and every $\lambda^j$. Let ${\cal O}^j=\{x(\xi)_j|\xi \in {\cal U}\}$. and assume  that  $d(x, \cdot)$ is  quasi-convex  in  its
second argument for all fixed $x\in X$.   If   there  exist $x(\xi^1)_j, \dots x(\xi^m)_j \in {\cal O}^j$ such that
${\cal O}^j \subseteq conv\{x(\xi^1)_j \dots x(\xi^m)_j\}$
then the center $x^\star$ of $x(\xi^1)_j, \dots x(\xi^m)_j$ is an optimal solution which minimizes the distance to the points in the set ${\cal O}^j$.
\end{thm}

\medskip

\begin{thm}
Let $(P(\xi:\xi \in {\cal U}))$ be an uncertain multi-objective optimization problem with problem instances $P(\xi,\lambda^j)$ which
admit  a  unique  optimal  solution for     every     $\xi \in {\cal U}$ and every $\lambda^j$.     Let ${\cal U}$ be a bounded polyhedral set (a polytope), i.e., ${\cal U} = conv\{\xi^1,\xi^2, \dots \xi^m\}$.  Assume  that  $d(x, \cdot)$ is  quasi-convex  in  its
second argument for all fixed $x\in X$ , and that $x:{\cal U} \rightarrow X$ is affine linear. Then the center of $x(\xi^1)_j, \dots x(\xi^m)_j$ with respect to the recovery costs $d$ solves the problem in Step \ref{step:x*} in the case of ${\cal O}^j=\{x(\xi)_j|\xi \in {\cal U}\}$.
\end{thm}

\medskip

\noindent
{\bf Remark.} The previous theorem describes  a case in which the extreme points
of the set ${ \cal O}^j$ can be determined without explicitly computing ${ \cal O}^j$,
but just by looking at set ${\cal U}$. As we can see, these extreme points do not rely on the chosen $\lambda^j$ value, yet they are determined from  set ${\cal U}$. Hence, we may conclude that this theorem shows how to choose a subset of scenarios in Step \ref{step:scenarios}, so that we end up with an accurate approximation of the initial uncertainty set when certain conditions are fulfilled. Thus, we finish with the following corollary.

\medskip

\begin{cor}
Let $(P(\xi:\xi \in {\cal U}))$ be an uncertain multi-objective optimization problem with problem instances $P(\xi,\lambda^j)$ which
admit  a  unique  optimal  solution for     every     $\xi \in {\cal U}$ and every $\lambda^j$.     Let ${\cal U}$ be a bounded polyhedral set (a polytope), i.e. ${\cal U} = conv\{\xi^1,\xi^2, \dots \xi^m\}$.  Assume  that  $d(x, \cdot)$ is  quasi-convex  in  its
second argument for all fixed $x\in X$ , and that $x:{\cal U} \rightarrow X$ is affine linear. Then, without loss of generality a set $\{\xi^1,\xi^2, \dots \xi^m \}$ may be chosen in Step \ref{step:scenarios} instead of the set ${\cal U}$.
\end{cor}

\medskip

Notice that, analogously to Algorithm \ref{alg:robustset}, the algorithm for finding recovery-to-efficiency robust  solutions by solving Problem $(\text{Rec-Eff center\_opt})$ (resp. $(\text{Rec-Eff median\_opt})$) may be deduced with the necessary changes. Indeed, it is enough to use $(\text{Rec-Eff center\_opt})$ (resp. $(\text{Rec-Eff median\_opt})$) in order to determine the set $X$ in Step~\ref{step:x*}.

\section{Experimental Analysis}
\label{sec:exp}

In this section, we are interested in comparing the recovery-to-efficiency robust sets obtained within Algorithm~\ref{alg:robustset} for 
problems $(\text{Rec-Eff center})$, $(\text{Rec-Eff median}$, $(\text{Rec-Eff center\_opt})$, and $(\text{Rec-Eff median\_opt})$. For testing purposes, we consider the following bi-objective knapsack problem: 
%
\begin{align}\label{prob}
\max~			&	\Big( \sum_{j=1}^\ell c^1_jx_j, \sum_{j=1}^\ell c^2_jx_j  \Big)		\tag{KP}	\\
\text{subject to~}	&	\sum_{j=1}^\ell w_jx_j\leq W		\notag	\\
				&	x_j \in x \in\{0,1\}^\ell				\notag
\end{align}
%
%
We consider $60$ instances of the bi-objective knapsack problem with different sizes. The number of items in the instances varies from $50$ to $500$. For each setting, $10$ different instances are independently generated. The costs of each item $c^1_j$  and $c^2_j$ as well as the weights $w_j$ are generated as random integer numbers
in $\llbracket 10,100 \rrbracket$. We assume that the objective coefficients $c^1_j$  and $c^2_j$  are subject to uncertainty, and thus 10 different realizations of these values are randomly generated (each realization corresponds to one scenario). The capacity of the knapsack is set as $W=\lceil \sum_{j=1}^\ell w_j/2 \rceil$ in all scenarios.  In all settings, the set $\Lambda$ is defined as $\Lambda=\{(0.00 + 0.01 \cdot k, 1 - 0.01 \cdot k) ~\vert~ k \in \set{0,1,\dots,100}\}$. As a solver we use the CPLEX 12.9. MIP solver which is able to optimally solve single objective knapsack instances in short time (less than 60 seconds). The distance between solution pairs is taken as the Hamming distance between binary strings.
To generate efficient solutions, we use weighted sum and Chebyshev scalarizing techniques. In the later case, we use the reference points $h_1$ and $h_2$, defined as

$$h_i=\max_ { \xi \in {\cal U}} \{\max c^i(\xi)_j x_j : \sum_{j=1}^\ell w_jx_j\leq W, x \in\{0,1\} ^\ell \}, \, , i=1,2$$
where ${\cal U}$ represents set of 10 different scenarios, and $c^i(\xi)_j, \, i=\set{1,2}, j=\set{1,2,\dots, \ell}$ realization of objective function values in the scenario $ \xi \in {\cal U}$.

\subsection{Comparing Robust Sets}

\begin{figure}[!t]
$$
\begin{matrix}
  \includegraphics [width = 7.5cm, angle = 0]{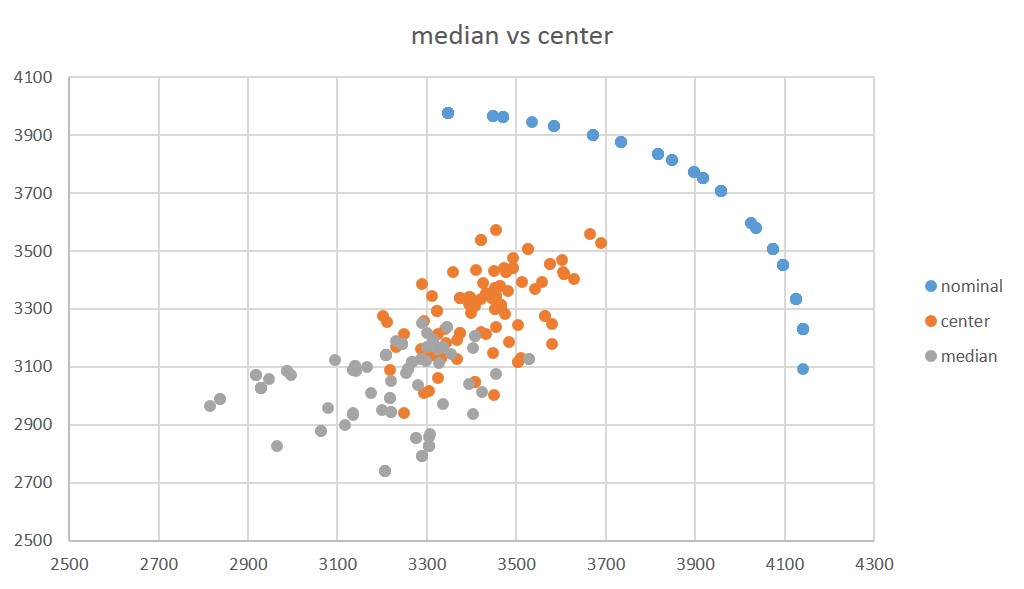} &
  \includegraphics [width = 7.5cm, angle = 0]{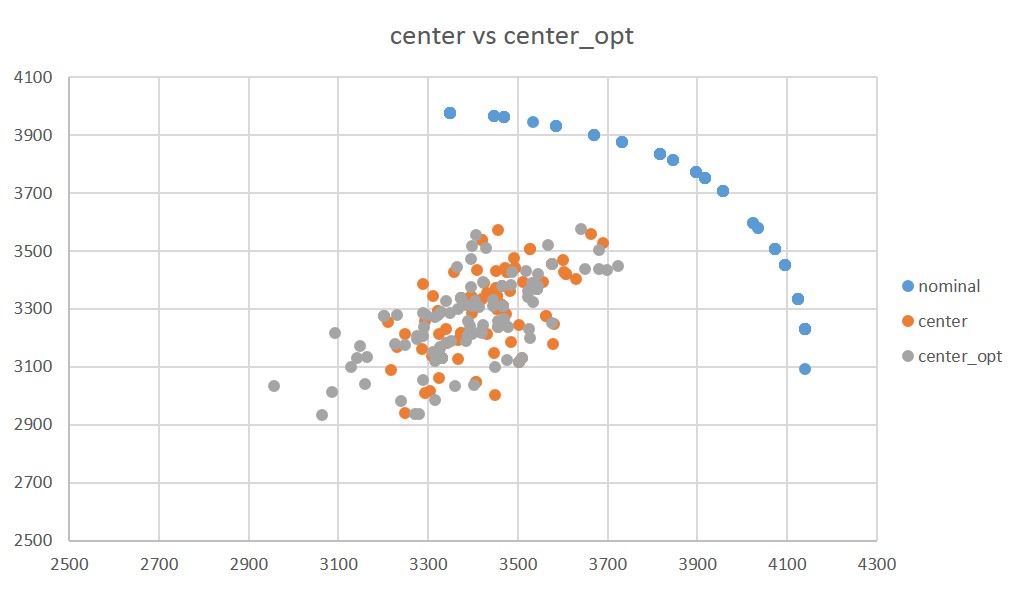} \\
  \includegraphics [width = 7.5cm, angle = 0]{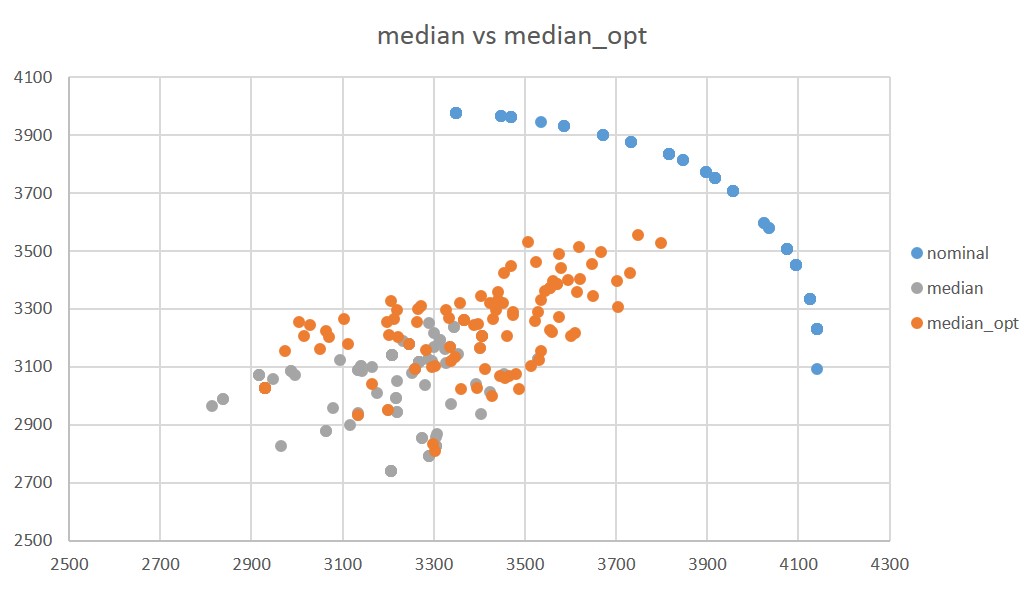} &\\
  \end{matrix}
  $$
\caption{Comparison of different robust sets for the weighted sum scalarizing function. \textsf{nominal} denotes the points from the efficient set in the nominal case; \textsf{center}, \textsf{center\_opt}, \textsf{median}, and \textsf{median\_opt} denote the set of recovery-to-efficiency robust solutions obtained considering problems
$(\text{Rec-Eff center})$,  $(\text{Rec-Eff center\_opt})$, $(\text{Rec-Eff median})$,  and $(\text{Rec-Eff median\_opt})$, respectively.} \label{fig:1}
\end{figure}

In the first set of experiments, we are interested in appreciating the differences among the efficient set of the nominal case (i.e., the first scenario) and the robust sets obtained by Algorithm \ref{alg:robustset} and the different models presented in Section~\ref{sec:rte}. 
They are presented in Fig.~\ref{fig:1}.
Each point in the graph represents the objective vector of each solution in the corresponding nominal case. More precisely, the solution $x$ from each set is represented as a point $(\sum_{j=1}^\ell c^1_jx_j, \sum_{j=1}^\ell c^2_jx_j)$ where $c^1_j$ and $c^2_j$ represent the objective coefficients in the nominal case. For testing purposes, we report results for the first instance of size $\ell=100$ from the previously described benchmark set.

\begin{figure}[!t]
$$
\begin{matrix}
  \includegraphics [width = 7.5cm, angle = 0]{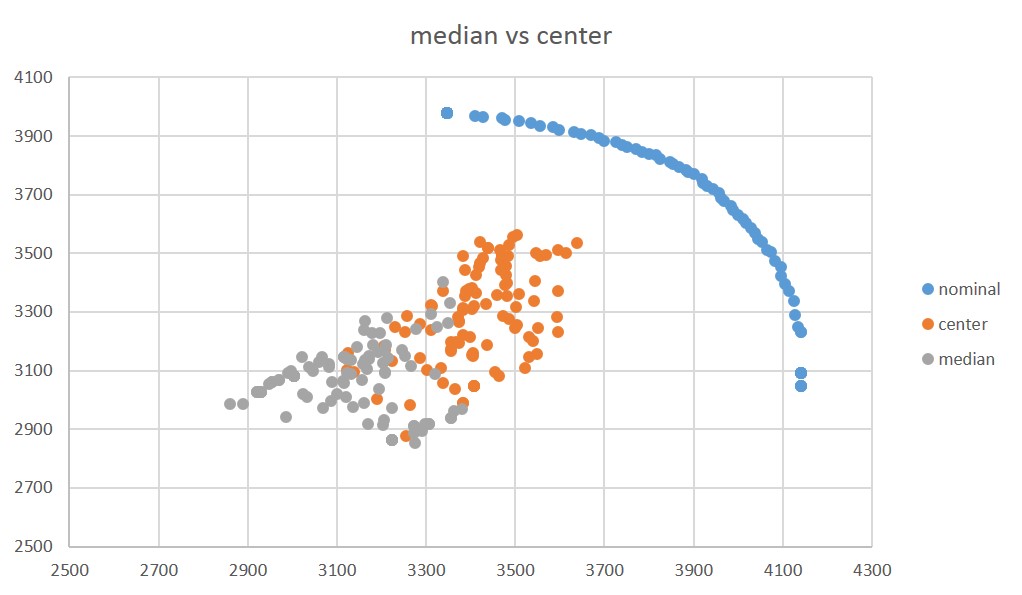} &
  \includegraphics [width = 7.5cm, angle = 0]{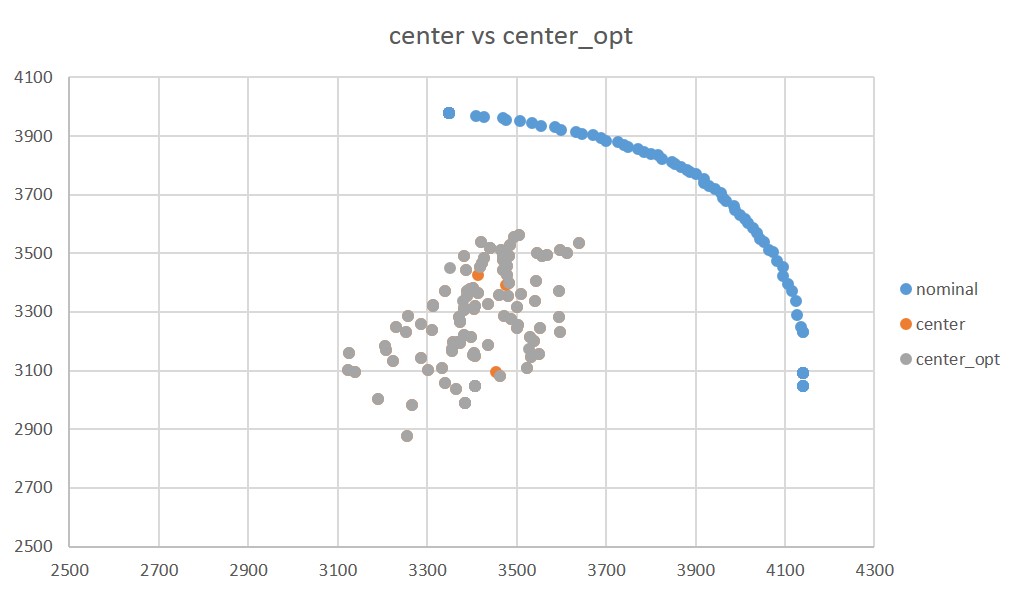} \\
  \includegraphics [width = 7.5cm, angle = 0]{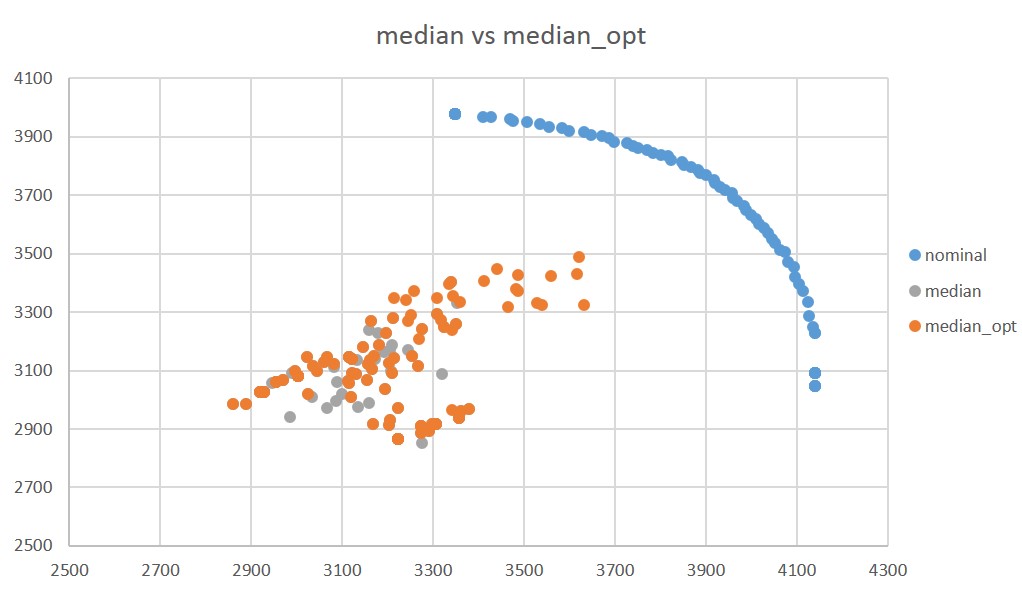} &\\
  \end{matrix}
  $$
\caption{Comparison of different robust sets for the Chebyshev scalarizing function. \textsf{nominal} denotes the points from the efficient set in the nominal case; \textsf{center}, \textsf{center\_opt}, \textsf{median}, and \textsf{median\_opt} denote the set of recovery-to-efficiency robust solutions obtained considering problems
$(\text{Rec-Eff center})$,  $(\text{Rec-Eff center\_opt})$, $(\text{Rec-Eff median})$,  and $(\text{Rec-Eff median\_opt})$, respectively.} \label{fig:2}
  \end{figure}

From Figs.~\ref{fig:1} and~\ref{fig:2}, several interesting observations may be derived. First, as expected, the robust sets do not follow the shape of the Pareto font.
Indeed, 
robust solutions are obtained as the closest ones to the solutions belonging to the efficient set for different scenarios. Comparing the sets obtained when solving $(\text{Rec-Eff center})$ and $(\text{Rec-Eff median})$ in Fig.~\ref{fig:1}--top, 
we observe that the points obtained for $(\text{Rec-Eff center})$ are closer to the nominal Pareto front than the ones obtained for $(\text{Rec-Eff median})$  in this particular case. A possible explanation for such outcome is that, in $(\text{Rec-Eff center})$, we optimize the worst distance. As a consequence, points tend to get closer to the nominal Pareto front.
By comparing the sets obtained for $(\text{Rec-Eff center} (\lambda))$ and $(\text{Rec-Eff center\_opt}(\lambda))$ in Fig.~\ref{fig:1}--middle, 
it follows that these two sets are different, although the objective vectors from both problems coincide for each~$\lambda$.
Such outcome may be justified by the fact that, in min-max types of problems such as these two, there are usually many solutions that produce the same optimal value.
Finally, by comparing the sets obtained for $(\text{Rec-Eff median} (\lambda))$ and $(\text{Rec-Eff median\_opt}(\lambda))$ in Fig.~\ref{fig:1}--bottom,
we see that the set generated solving $(\text{Rec-Eff median\_opt}(\lambda))$ tend to get closer to the nominal Pareto front. 
However, we will see in the next section that
solving $(\text{Rec-Eff median\_opt}(\lambda))$ rather than $(\text{Rec-Eff median}(\lambda))$ actually leads to cheaper total recovery~costs in all test cases.

\subsection{Comparing $(\text{Rec-Eff median})$ and $(\text{Rec-Eff median\_opt})$}
In this section, we compare the total recovery cost induced when solving $(\text{Rec-Eff median})$ and $(\text{Rec-Eff median\_opt})$. In particular, we are interested in appreciating to which extent $(\text{Rec-Eff median\_opt})$ actually yields to a smaller total recovery cost than $(\text{Rec-Eff median})$. A summary of our experimental results is presented in
Table \ref{tabFB1} (weighted sum scalarizing function used to generate Pareto efficient solutions) and Table \ref{tabFB3} (Chebyshev scalarizing function used to generate Pareto efficient solutions).
Columns `$\ell$' and `$\#$' respectively report the number of items in the considered instances and the number of tested instances. Columns `median' and `median\_opt' respectively report the total recovery costs obtained considering $(\text{Rec-Eff median})$ and $(\text{Rec-Eff median\_opt})$. Column `diff' reports the difference between values reported  under columns `median' and `median\_opt', while Column `\% dev.' reports the percentage deviation of the `median' value from the `median\_opt' value.


\begin{table}[p]
\caption{Median vs Median\_opt: weighted sum scalarizing function. } \label{tabFB1}
\begin{center}
\begin{small}
\begin{tabular}{cc|rr|r||cc|rr|r}

 & & \multicolumn{2}{c|}{recovery cost} & & & & \multicolumn{2}{c|}{recovery cost} & \\
 \hline
$\ell$ & \# & median & median\_opt & \% dev. & $\ell$ & \# & median & median\_opt& \% dev.\\
\hline
50	&	1	&	11016	&	11008	&	0.07	&	300	&	1	&	65785	&	65751	&	0.05	\\
50	&	2	&	10898	&	10896	&	0.02	&	300	&	2	&	62926	&	62904	&	0.03	\\
50	&	3	&	12152	&	12146	&	0.05	&	300	&	3	&	63154	&	63091	&	0.10	\\
50	&	4	&	10610	&	10602	&	0.08	&	300	&	4	&	60785	&	60748	&	0.06	\\
50	&	5	&	10455	&	10450	&	0.05	&	300	&	5	&	64223	&	64188	&	0.05	\\
50	&	6	&	11926	&	11923	&	0.03	&	300	&	6	&	64852	&	64838	&	0.02	\\
50	&	7	&	10913	&	10910	&	0.03	&	300	&	7	&	57741	&	57707	&	0.06	\\
50	&	8	&	11155	&	11151	&	0.04	&	300	&	8	&	63348	&	63304	&	0.07	\\
50	&	9	&	11298	&	11286	&	0.11	&	300	&	9	&	60608	&	60554	&	0.09	\\
50	&	10	&	11024	&	11011	&	0.12	&	300	&	10	&	67405	&	67367	&	0.06	\\
\hline
\multicolumn{2}{c|}{Average:}&	11144.7	&	11138.3	&	0.06	&	\multicolumn{2}{c|}{Average:} &	63082.7	&	63045.2	&		0.06	\\
\hline
100	&	1	&	24404	&	24399	&	0.02	&	400	&	1	&	82004	&	81961	&	0.05	\\
100	&	2	&	23454	&	23432	&	0.09	&	400	&	2	&	87684	&	87633	&	0.06	\\
100	&	3	&	21571	&	21551	&	0.09	&	400	&	3	&	89224	&	89165	&	0.07	\\
100	&	4	&	23335	&	23318	&	0.07	&	400	&	4	&	85541	&	85495	&	0.05	\\
100	&	5	&	22285	&	22274	&	0.05	&	400	&	5	&	86087	&	86031	&	0.07	\\
100	&	6	&	20110	&	20103	&	0.03	&	400	&	6	&	82055	&	81996	&	0.07	\\
100	&	7	&	18994	&	18981	&	0.07	&	400	&	7	&	83132	&	83082	&	0.06	\\
100	&	8	&	22310	&	22288	&	0.10	&	400	&	8	&	83291	&	83224	&	0.08	\\
100	&	9	&	23386	&	23366	&	0.09	&	400	&	9	&	88815	&	88758	&	0.06	\\
100	&	10	&	21195	&	21182	&	0.06	&	400	&	10	&	84633	&	84584	&	0.06	\\
\hline
\multicolumn{2}{c|}{Average:} &	22104.4	&	22089.4	&	0.07	&\multicolumn{2}{c|}{Average:} &	85246.6	&	85192.9	&	0.06 \\
\hline
200	&	1	&	44127	&	44116	&	0.02	&	500	&	1	&	111054	&	111000	&	0.05	\\
200	&	2	&	40290	&	40250	&	0.10	&	500	&	2	&	109701	&	109639	&	0.06	\\
200	&	3	&	42938	&	42908	&	0.07	&	500	&	3	&	115990	&	115930	&	0.05	\\
200	&	4	&	41546	&	41519	&	0.07	&	500	&	4	&	108417	&	108363	&	0.05	\\
200	&	5	&	45840	&	45815	&	0.05	&	500	&	5	&	108518	&	108454	&	0.06	\\
200	&	6	&	45029	&	45007	&	0.05	&	500	&	6	&	106668	&	106598	&	0.07	\\
200	&	7	&	43516	&	43464	&	0.12	&	500	&	7	&	111305	&	111238	&	0.06	\\
200	&	8	&	43579	&	43553	&	0.06	&	500	&	8	&	106073	&	106010	&	0.06	\\
200	&	9	&	44524	&	44506	&	0.04	&	500	&	9	&	107972	&	107917	&	0.05	\\
200	&	10	&	39809	&	39772	&	0.09	&	500	&	10	&	108015	&	107950	&	0.06	\\
\hline
\multicolumn{2}{c|}{Average:} &	43119.8	&	43091.0	&	0.07	&\multicolumn{2}{c|}{Average:} &	109371.3	&	109309.9	&	0.06	\\
\hline
\end{tabular}
\end{small}
\end{center}
\end{table}

\begin{table}[p]
\caption{Median vs Median\_opt: Chebyshev scalarizing function. } \label{tabFB3}
\begin{center}
\begin{small}
\begin{tabular}{cc|rr|r||cc|rr|r}
 & & \multicolumn{2}{c|}{recovery cost} & & & & \multicolumn{2}{c|}{recovery cost} & \\
 \hline
$\ell$ & \# & median & median\_opt& \% dev. & $\ell$ & \# & median & median\_opt & \% dev.\\
\hline
50	&	1	&	11913	&	11785	&	1.09	&	300	&	1	&	66537	&	66151	&	0.58	\\
50	&	2	&	11824	&	11796	&	0.24	&	300	&	2	&	63720	&	63421	&	0.47	\\
50	&	3	&	12823	&	12756	&	0.53	&	300	&	3	&	64661	&	64119	&	0.85	\\
50	&	4	&	11919	&	11766	&	1.30	&	300	&	4	&	61821	&	61430	&	0.64	\\
50	&	5	&	11256	&	11202	&	0.48	&	300	&	5	&	65976	&	65657	&	0.49	\\
50	&	6	&	12766	&	12578	&	1.49	&	300	&	6	&	65523	&	65239	&	0.44	\\
50	&	7	&	11495	&	11465	&	0.26	&	300	&	7	&	58860	&	58553	&	0.52	\\
50	&	8	&	12024	&	11986	&	0.32	&	300	&	8	&	63301	&	62977	&	0.51	\\
50	&	9	&	12266	&	12233	&	0.27	&	300	&	9	&	63362	&	62899	&	0.74	\\
50	&	10	&	11749	&	11709	&	0.34	&	300	&	10	&	68157	&	67715	&	0.65	\\
\hline
\multicolumn{2}{c|}{Average:}&	12003.5	&	11927.6	&	0.63	&	\multicolumn{2}{c|}{Average:} &64191.8	&	63816.1	&	0.59	\\
\hline
100	&	1	&	25392	&	25343	&	0.19	&	400	&	1	&	81841	&	81340	&	0.62	\\
100	&	2	&	24701	&	24586	&	0.47	&	400	&	2	&	88965	&	88385 &	0.66	\\
100	&	3	&	22762	&	22676	&	0.38	&	400	&	3	&	89201	&	88774	&	0.48	\\
100	&	4	&	24232	&	24026	&	0.86	&	400	&	4	&	87971	&	87359	&	0.70	\\
100	&	5	&	23370	&	23279	&	0.39	&	400	&	5	&	89514	&	89114	&	0.45	\\
100	&	6	&	22577	&	22524	&	0.24	&	400	&	6	&	84051	&	83484	&	0.68	\\
100	&	7	&	19748	&	19678	&	0.36	&	400	&	7	&	84872	&	84430	&	0.52	\\
100	&	8	&	23410	&	23298	&	0.48	&	400	&	8	&	88335	&	88032	&	0.34	\\
100	&	9	&	24646	&	24415	&	0.95	&	400	&	9	&	86183	&	85952	&	0.27	\\
100	&	10	&	22767	&	22528	&	1.06	&	400	&	10	&	91093	&	90626	&	0.52	\\
\hline
\multicolumn{2}{c|}{Average:} &	23360.5	&	23235.3	&	0.54	&\multicolumn{2}{c|}{Average:} &	87202.6	&	86749.6	&	0.52\\
\hline
200	&	1	&	45004	&	44882	&	0.27	&	500	&	1	&	113140	&	112902	&	0.21	\\
200	&	2	&	41356	&	41117	&	0.58	&	500	&	2	&	106496	&	106284	&	0.20	\\
200	&	3	&	44639	&	44280	&	0.81	&	500	&	3	&	112648	&	112230	&	0.37	\\
200	&	4	&	43266	&	42848	&	0.98	&	500	&	4	&	110058	&	109581	&	0.44	\\
200	&	5	&	45828	&	45611	&	0.48	&	500	&	5	&	107168	&	107029	&	0.13	\\
200	&	6	&	45934	&	45706	&	0.50	&	500	&	6	&	109258	&	109004	&	0.23	\\
200	&	7	&	44294	&	44015	&	0.63	&	500	&	7	&	115781	&	115652	&	0.11	\\
200	&	8	&	44771	&	44473	&	0.67	&	500	&	8	&	108145	&	107913	&	0.21	\\
200	&	9	&	45630	&	45402	&	0.50	&	500	&	9	&	101011	&	100807	&	0.20	\\
200	&	10	&	41135	&	40830	&	0.75	&	500	&	10	&	109448	&	109202	&	0.23	\\

\hline
\multicolumn{2}{c|}{Average:} &	44185.7	&	43916.4	&	0.62	&\multicolumn{2}{c|}{Average:} &	109315.3	&	109060.4	&	0.23	\\
\hline
\end{tabular}
\end{small}
\end{center}
\end{table}

From the reported results, it follows that considering $(\text{Rec-Eff median\_opt})$ instead of $(\text{Rec-Eff median})$ yield certain savings on each test instance regardless of the scalarizing techniques used. In the case a weighted sum scalarizing function is used, such one saving may be up to 0.12\%, while in the case of Chebyshev, the savings may go up to 1.49\%. We also observe, that average savings  over instances of the same size in the case weighted sum scalarization is used are smaller than in the case  Chebyshev scalarization is used. In particular, if the weighted sum scalarizing function is used, the average savings  are 0.06\% except on instances with $\ell=100$ and $\ell=200$ where average savings are 0.07\%, while if a Chebyshev scalarizing function is used the average savings are at least 0.23\%. So, we may conclude that $(\text{Rec-Eff median\_opt})$ is better option than $(\text{Rec-Eff median})$ in the cases where problems $P(\xi,\lambda)$ have more than one optimal solution as it is case in this paper. In addition, it follows that weighted sum has less impact on the difference of objective function values of $(\text{Rec-Eff median\_opt})$ and $(\text{Rec-Eff median})$, than Chebyshev.


\subsection{Comparing $(\text{Rec-Eff center})$ and $(\text{Rec-Eff center\_opt})$}
In this section, we compare the total recovery cost induced when solving problems $(\text{Rec-Eff center})$ and $(\text{Rec-Eff center\_opt})$ in order  to determine savings yielded by using $(\text{Rec-Eff center\_opt})$ instead of $(\text{Rec-Eff center})$. A summary of our experimental results is presented in
Table \ref{tabFB2} (weighted scalarizing function used to generate Pareto efficient solutions) and Table \ref{tabFB4}  (Chebyshev scalarizing function used to generate Pareto efficient solutions).
Columns `center' and `center\_opt' respectively report the total recovery cost obtained considering $(\text{Rec-Eff center})$ and $(\text{Rec-Eff center\_opt}$.

From the results reported in Table \ref{tabFB2}, it follows that on 6 out of 60 instances total recovery costs produced considering $(\text{Rec-Eff center\_opt})$ and $(\text{Rec-Eff center})$ are the same. On smallest instances with $\ell=50$, the total recovery costs coincide on 3 out of 10 cases. On the other hand the savings induced considering $(\text{Rec-Eff center\_opt})$ instead of $(\text{Rec-Eff center})$ may be up to 0.23\%. In the case that a Chebyshev scalarizing function is used to generate Pareto efficient solutions (Table \ref{tabFB4}), there is no instance on which total recovery costs produced considering $(\text{Rec-Eff center\_opt})$ and $(\text{Rec-Eff center})$ are the same. In addition, we observe that  considering $(\text{Rec-Eff center\_opt})$ instead of $(\text{Rec-Eff center})$ may yield savings of up to 1.99\%. So, we may conclude that $(\text{Rec-Eff center\_opt})$ is better option than $(\text{Rec-Eff center})$ either weighted sum or Chebyshev is used to generate Pareto efficient solutions. In addition, it follows that the weighted sum scalarizing function has less impact on the difference of objective function values of $(\text{Rec-Eff cnter\_opt})$ and $(\text{Rec-Eff center})$, than Chebyshev.

\begin{table}[p]
\caption{Center vs center\_opt: weighted sum scalarizing function.} \label{tabFB2}
\begin{center}
\begin{small}
\begin{tabular}{cc|rr|r||cc|rr|r}
$\ell$& \# &center & center\_opt& \% dev. & $\ell$ & \# & center & center\_opt& \% dev.\\
\hline
50	&	1	&	1308	&	1307	&	0.08	&	300	&	1	&	6780	&	6780	&	0.00	\\
50	&	2	&	1275	&	1275	&	0.00	&	300	&	2	&	6542	&	6540	&	0.03	\\
50	&	3	&	1406	&	1405	&	0.07	&	300	&	3	&	6552	&	6547	&	0.08	\\
50	&	4	&	1280	&	1279	&	0.08	&	300	&	4	&	6312	&	6311	&	0.02	\\
50	&	5	&	1253	&	1253	&	0.00	&	300	&	5	&	6569	&	6566	&	0.05	\\
50	&	6	&	1380	&	1379	&	0.07	&	300	&	6	&	6649	&	6648	&	0.02	\\
50	&	7	&	1342	&	1342	&	0.00	&	300	&	7	&	5956	&	5953	&	0.05	\\
50	&	8	&	1337	&	1336	&	0.07	&	300	&	8	&	6503	&	6499	&	0.06	\\
50	&	9	&	1331	&	1328	&	0.23	&	300	&	9	&	6240	&	6238	&	0.03	\\
50	&	10	&	1298	&	1296	&	0.15	&	300	&	10	&	6916	&	6912	&	0.06	\\
\hline
\multicolumn{2}{c|}{Average:} 	&	1321.0	&	1320.0	&	0.08	&\multicolumn{2}{c|}{Average:} &	6501.9	&	6499.4	&	0.04	\\
\hline
100	&	1	&	2612	&	2612	&		0.00	&	400	&	1	&	8356	&	8354	&	0.02	\\
100	&	2	&	2519	&	2516	&		0.12	&	400	&	2	&	8894	&	8890	&	0.04	\\
100	&	3	&	2364	&	2362	&		0.08	&	400	&	3	&	9063	&	9059	&	0.04	\\
100	&	4	&	2582	&	2581	&	    0.04	&	400	&	4	&	8726	&	8723	&	0.03	\\
100	&	5	&	2461	&	2459	&		0.08	&	400	&	5	&	8833	&	8832	&	0.01	\\
100	&	6	&	2250	&	2249	&		0.04	&	400	&	6	&	8366	&	8362	&	0.05	\\
100	&	7	&	2066	&	2065	&		0.05	&	400	&	7	&	8481	&	8478	&	0.04	\\
100	&	8	&	2438	&	2436	&		0.08	&	400	&	8	&	8497	&	8494	&	0.04	\\
100	&	9	&	2507	&	2506	&		0.04	&	400	&	9	&	9095	&	9092	&	0.03	\\
100	&	10	&	2326	&	2324	&		0.09	&	400	&	10	&	8610	&	8606	&	0.05	\\
\hline
\multicolumn{2}{c|}{Average:} 	&	2412.5	&	2411.0	&	0.06	&		\multicolumn{2}{c|}{Average:} &	8692.1	&	8689.0	&	0.04	\\
\hline
200	&	1	&	4625	&	4624	&	0.02	&	400	&	10	&	8610	&	8606	&	0.05	\\
200	&	2	&	4205	&	4204	&	0.02	&	500	&	1	&	11285	&	11278	&	0.06	\\
200	&	3	&	4529	&	4527	&	0.04	&	500	&	2	&	11143	&	11140	&	0.03	\\
200	&	4	&	4435	&	4432	&	0.07	&	500	&	3	&	11766	&	11766	&	0.00	\\
200	&	5	&	4762	&	4759	&	0.06	&	500	&	4	&	11024	&	11023	&	0.01	\\
200	&	6	&	4644	&	4642	&	0.04	&	500	&	5	&	10994	&	10989	&	0.05	\\
200	&	7	&	4553	&	4549	&	0.09	&	500	&	6	&	10792	&	10789	&	0.03	\\
200	&	8	&	4572	&	4569	&	0.07	&	500	&	7	&	11260	&	11254	&	0.05	\\
200	&	9	&	4623	&	4621	&	0.04	&	500	&	8	&	10801	&	10798	&	0.03	\\
200	&	10	&	4219	&	4216	&	0.07	&	500	&	9	&	10942	&	10937	&	0.05	\\
\hline
\multicolumn{2}{c|}{Average:} &	4516.7	&	4514.3	&	0.05	&		\multicolumn{2}{c|}{Average:} 	&	11091.3	&	11087.8	&	0.03	\\

\hline
\end{tabular}
\end{small}
\end{center}
\end{table}



\begin{table}[p]
\caption{Center vs center\_opt: Chebyshev scalarizing function.} \label{tabFB4}
\begin{center}
\begin{small}
\begin{tabular}{cc|rr|r||cc|rr|r}
$\ell$ & \# &center & center\_opt&  \% dev. & $\ell$ & \# & center & center\_opt& \% dev.\\
\hline
50	&	1	&	1422	&	1403	&	1.35	&	300	&	1	&	6965	&	6942	&	0.33	\\
50	&	2	&	1380	&	1373	&	0.51	&	300	&	2	&	6745	&	6740	&	0.07	\\
50	&	3	&	1462	&	1460	&	0.14	&	300	&	3	&	6798	&	6773	&	0.37	\\
50	&	4	&	1445	&	1430	&	1.05	&	300	&	4	&	6501	&	6486	&	0.23	\\
50	&	5	&	1329	&	1323	&	0.45	&	300	&	5	&	6838	&	6823	&	0.22	\\
50	&	6	&	1487	&	1458	&	1.99	&	300	&	6	&	6815	&	6768	&	0.69	\\
50	&	7	&	1398	&	1395	&	0.22	&	300	&	7	&	6151	&	6138	&	0.21	\\
50	&	8	&	1409	&	1406	&	0.21	&	300	&	8	&	6555	&	6537	&	0.28	\\
50	&	9	&	1456	&	1455	&	0.07	&	300	&	9	&	6591	&	6561	&	0.46	\\
50	&	10	&	1410	&	1405	&	0.36	&	300	&	10	&	7099	&	7063	&	0.51	\\
\hline
\multicolumn{2}{c|}{Average:} 	&	1419.8	&	1410.8	&	0.63	&\multicolumn{2}{c|}{Average:} &	6705.8	&	6683.1	&	0.34	\\
\hline
100	&	1	&	2726	&	2723	&	0.11	&	400	&	1	&	8431	&	8412	&	0.23	\\
100	&	2	&	2662	&	2654	&	0.30	&	400	&	2	&	9103	&	9083	&	0.22	\\
100	&	3	&	2526	&	2520	&	0.24	&	400	&	3	&	9150	&	9137	&	0.14	\\
100	&	4	&	2698	&	2695	&	0.11	&	400	&	4	&	9067	&	9041	&	0.29	\\
100	&	5	&	2606	&	2598	&	0.31	&	400	&	5	&	9254	&	9248	&	0.06	\\
100	&	6	&	2495	&	2495	&	0.00	&	400	&	6	&	8663	&	8638	&	0.29	\\
100	&	7	&	2194	&	2191	&	0.14	&	400	&	7	&	8683	&	8671	&	0.14	\\
100	&	8	&	2578	&	2574	&	0.16	&	400	&	8	&	9102	&	9070	&	0.35	\\
100	&	9	&	2624	&	2620	&	0.15	&	400	&	9	&	8894	&	8876	&	0.20	\\
100	&	10	&	2515	&	2490	&	1.00	&	400	&	10	&	9342	&	9317	&	0.27	\\
\hline
\multicolumn{2}{c|}{Average:} 	&	2562.4	&	2556	&	0.25	&		\multicolumn{2}{c|}{Average:} &	8968.9	&	8949.3	&	0.22	\\
\hline
200	&	1	&	4751	&	4749	&	0.04	&	500	&	1	&	11588	&	11579	&	0.08	\\
200	&	2	&	4399	&	4391	&	0.18	&	500	&	2	&	11167	&	11135	&	0.29	\\
200	&	3	&	4753	&	4739	&	0.30	&	500	&	3	&	11819	&	11807	&	0.10	\\
200	&	4	&	4664	&	4627	&	0.80	&	500	&	4	&	11040	&	11026	&	0.13	\\
200	&	5	&	4815	&	4809	&	0.12	&	500	&	5	&	10882	&	10868	&	0.13	\\
200	&	6	&	4751	&	4746	&	0.11	&	500	&	6	&	11525	&	11496	&	0.25	\\
200	&	7	&	4664	&	4645	&	0.41	&	500	&	7	&	11229	&	11199	&	0.27	\\
200	&	8	&	4774	&	4758	&	0.34	&	500	&	8	&	10904	&	10893	&	0.10	\\
200	&	9	&	4720	&	4714	&	0.13	&	500	&	9	&	11287	&	11268	&	0.17	\\
200	&	10	&	4456	&	4437	&	0.43	&	500	&	10	&	10350	&	10336	&	0.14	\\

\hline
\multicolumn{2}{c|}{Average:} &	4674.7	&	4661.5	&	0.29	&		\multicolumn{2}{c|}{Average:} 	&	11179.1	&	11160.7	&	0.16	\\

\hline
\end{tabular}
\end{small}
\end{center}
\end{table}

\subsection{Influence of $\epsilon$ tolerance}

\begin{figure}[H]
$$
\begin{matrix}
  \includegraphics [width = 8cm, angle = 0]{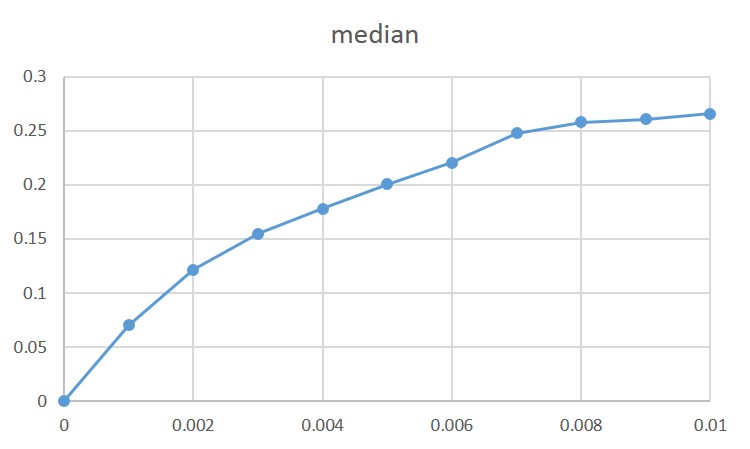} &
  \includegraphics [width = 8cm, angle = 0]{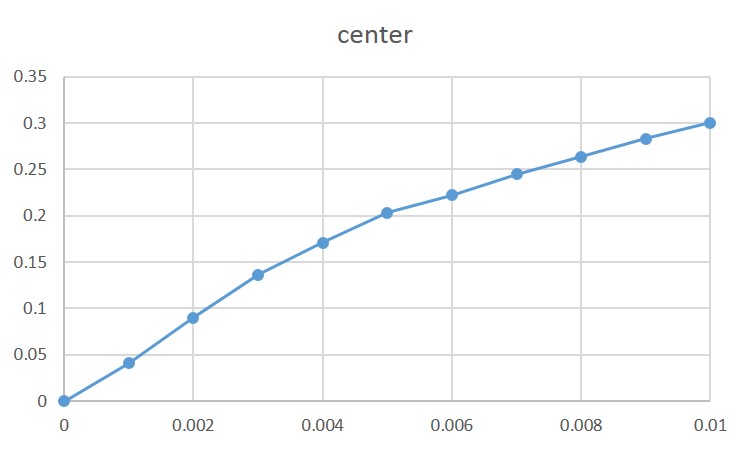} \\
  \end{matrix}
  $$
    \caption{Sensitivity to $\epsilon$: weighted sum scalarizing function.} \label{fig:3}
  \end{figure}

\begin{figure}[H]
$$
\begin{matrix}
  \includegraphics [width = 8cm, angle = 0]{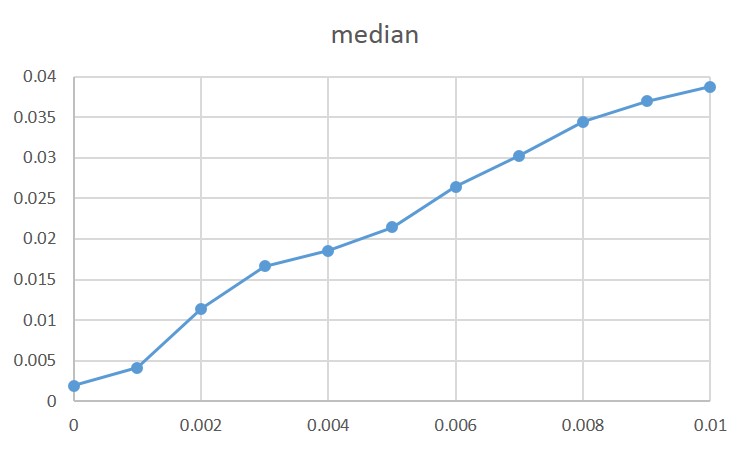} &
  \includegraphics [width = 8cm, angle = 0]{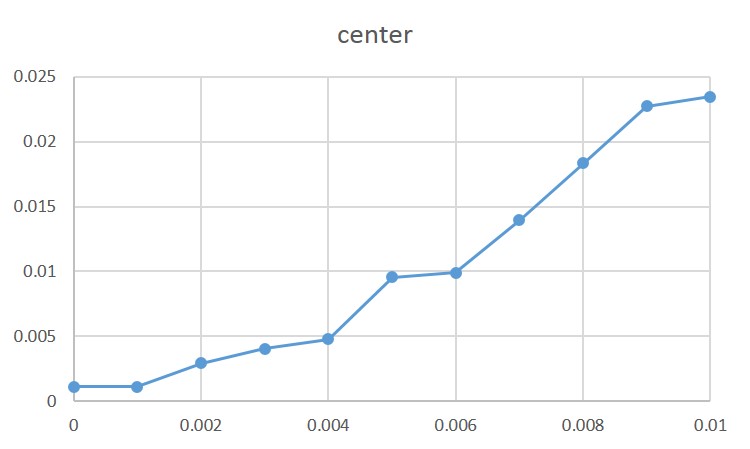} \\
  \end{matrix}
  $$
\caption{Sensitivity to $\epsilon$: Chebyshev scalarizing function.} \label{fig:4}
  \end{figure}

In this section, we consider the problems $(\text{Rec-Eff center\_opt})$ and $(\text{Rec-Eff median\_opt})$ where instead using the standard constraint
\begin{equation}\label{eq0}
g(x(\xi, \lambda))\leq g_{opt}(\xi, \lambda), \;  \; \xi \in {\cal U}
\end{equation}
we use constraint
\begin{equation}\label{eq1}
g(x(\xi, \lambda)) \leq (1+\epsilon) g_{opt}(\xi, \lambda), \;  \; \xi \in {\cal U}
\end{equation}
which allows that "efficient solutions" in scenarios are within some tolerance $\epsilon$ with respect to the optimal objective value. We consider 11 different values of $\epsilon$ ranging from 0 to 0.01 with step 0.001. (i.e., allowing percentage deviation of up to 1\% from the optimal). Since we have maximization problems, we allow that the objective function values of efficient solutions are up to 1\% less than optimal ones. a For testing purposes, we used the first instance with $\ell = 100$ from our benchmark set. In figures \ref{fig:3} and \ref{fig:4} we present the results for both center and median objectives. In
Figure \ref{fig:3}, Pareto efficient solutions for each scenario are generated using the weighted sum scalarizing technique,
while in In Figure \ref{fig:4} the Chebyshev scalarizing technique is used to do so. At each graph, on $x$-axis  $\epsilon$  values are given, while on  $y$-axis deviations of the value $v(\epsilon)$ from the value $v(0)$ are given, where $v(\epsilon)$ stands for the objective value of the corresponding problem where constraint \eqref{eq1} is used, while $v(\epsilon)$ stands for the objective value of the corresponding problem where constraint \eqref{eq0} is used. Deviations are computed using the following formula:
$$ \frac{v(0)-v(\epsilon)}{v(0)}.$$

From figures \ref{fig:3} and \ref{fig:4}, we may infer that regardless of the scalarizing technique used to generate Pareto efficient solutions, problems $(\text{Rec-Eff center\_opt})$ and $(\text{Rec-Eff median\_opt})$ are sensitive to the imposed value of $\epsilon$. Also, we observe that in any case increase of $\epsilon$ value decreases $v(\epsilon)$ value.

\section{Concluding Remarks}
\label{sec:conclu}
In this paper, we presented a new robustness concept for uncertain multi-objective optimization problems.
The so-obtained recovery-to-efficiency approach is inspired by existing robustness concepts from uncertain single-objective optimization called recovery-to-optimality \citep{goerigk2014recovery}.
As shown in the paper, recovery-to-efficiency actually extends and generalizes recovery-to-optimality, as the latter is as a special case of the former. The paper also proposes a number of approaches for determining a set of recovery-to-efficiency robust solutions. An extensive experimental analysis allowed us to disclose differences among robust sets obtained using different concepts as well as to deduce some interesting observations.
First, as expected, the robust sets do not
follow the shape of the Pareto font. Second, the weighted sum scalarizing technique has less impact on the difference of objective function values of
concepts $(\text{Rec-Eff center\_opt})$ and $(\text{Rec-Eff center})$, than Chebyshev. The same holds considering concepts $(\text{Rec-Eff median\_opt})$ and $(\text{Rec-Eff median})$. Finally, regardless of the scalarizing technique used to generate efficient
solutions, concepts $(\text{Rec-Eff center\_opt})$  and $(\text{Rec-Eff median\_opt})$  are sensitive to the imposed value of the $\epsilon$~parameter.

Future works may include proposing new concepts to determine recovery-to-efficiency robust sets, as well as proposing other robustness concepts for uncertain multi-objective optimization problems. Since, in the paper, a weighted-sum scalarizing technique is used, it would be interesting to develop recovery-to-efficiency robustness concepts for other multi-objective optimization techniques to identify efficient solutions, including other scalarizing functions. Moreover, the future work will focus on finding an industrial partner and the use of manuscript's content and concept in a real world problem.

\section*{References}
\bibliographystyle{elsarticle-harv}
\bibliography{robust}

\end{document}